# A self-contained theory of truth

David Sikter

**ABSTRACT**

Tarski's undefinability theorem states that a formal system based on conventional *predicate logic* (PL) cannot talk about its own truth predicate. PL is, however, not the only formal language imaginable. In this paper, it will be shown that it is possible to construct a formal system, based not on PL but a more restrictive formal language, which is *self-contained* in the sense that within this system, we can both talk about and even define the system's own truth predicate. This hints at new ways of understanding Russell's, Gödel's, and Tarski's discoveries, and new ways of tackling the vicious circles that give rise to these problems. This new system presents an interesting alternative foundational mathematical framework.

## Introduction

It has long been known that using language to talk about language – and especially using language to talk about the notion of *truth* – tends to lead to contradictions arising from vicious circularity, like different versions of the infamous *liar paradox*. The full extent of these problems became apparent in the early 20[th] century with the discovery of Russell's Paradox, Gödel's incompleteness theorem, and Tarski's undefinability theorem – discoveries that led to what has been called *the foundational crisis of mathematics*.

This crisis remains essentially unresolved to this day – there is still no consensus on how the results of Russell, Gödel, and Tarski should be interpreted, or how the vicious circles that give rise to these problems should be tackled.

Many seem to interpret these results as demonstrating that the notion of truth is fundamentally undefinable and that it is impossible to devise a complete theory of truth. I contend, however, that this is a premature conclusion. In this paper, I am going to demonstrate that it is actually perfectly possible to construct a formal notion of truth, which is both sensible, has a great deal of expressive power, and allows us to both talk about and even define the system's own truth predicate, without contradiction.

# What is truth?

To get to the bottom of these problems, we need to get a bit philosophical, and we need to start from the very beginning. We have this intuitive notion that there is something we call *truth* – but what, exactly, *is* truth?

The notion of truth that we normally engage with is, of course, that of *informal* human languages, like English. However, as many have noted, these informal languages are extremely vague and imprecise – they allow us to formulate ambiguous, nonsensical, and paradoxical statements.

Because of problems like these, many have argued that we need some kind of *ideal language* – a precise, formal language with well-defined semantics and truth-values. This idea was popularized in the early 20th century by philosophers like Bertrand Russell [1] and Ludwig Wittgenstein [2], and it seems to me that this is indeed the only sensible conclusion. In fact, it seems to me that if a certain notion of truth can *not* be formalized, then something must be fundamentally *wrong* with it. Why? Because if a certain set of ideas cannot be formalized, this means that all attempts to make them *precise* must run into some kind of *contradictions*, and that is, I submit, also exactly what we mean when we say that a set of ideas is *incoherent*. An *unformalizable* notion of truth is thus, I contend, an *incoherent* notion of truth (by definition, even, I would argue).

I am thus, in this paper, going to assume that we want a notion of truth that can be formalized, i.e. captured in some kind of formal system – and that notions of truth that can *not* be formalized must be rejected.

Of course, this begs the question of what exactly we mean by a "formal system". I am going to adopt the most generic definition of this concept that I can imagine, which is that a formal system consists of four components, namely:

1) a *semantic schema*, which defines the semantics of the system,
2) a *model*, which contains basic, atomic facts about the World, encoded according to the principles defined by the semantic schema,
3) a *formal grammar*, which defines what substitute a valid statement within the system, and finally
4) an *evaluation function* that defines the truth-values of valid statements.

## An analogy

To illustrate the general idea behind this definition, consider the analogy of a database system. To describe some aspects of the World using a database system, we first need to define our system's semantics. That is: we need to decide *what* we want it to describe and *how*. Then, of course, we need a database populated with the right data. The database is our system's model of the World – it should store information about the World, encoded according to the principles we have decided in our semantic schema.

When we have a populated database, we can ask queries against it. For that, we use a formal query language, like SQL. If a certain statement is a valid SQL query or not can be determined by a simple algorithm – the SQL parser. This algorithm defines our system's formal grammar. If a statement is valid, we can go on and evaluate it. For that, we again use an algorithm – the SQL evaluator. It requires two pieces of information: the data in the database and a string defining the statement to evaluate. From these two pieces of information, the correct answer to the query can be computed algorithmically. This algorithm is our system's evaluation function.

## The information-theoretic view of language

The analogy of a database system is a very interesting and powerful one. It suggests that we can define the notion of a formal system using three core concepts: *data structures*, *algorithms*, and *semantics*. I like to call this *the information-theoretic view of language*.

By a "data structure", I here mean *some kind of well-defined structure that can be used to encode information*. The model of our system (i.e. the database in our analogy) is, when it all comes down to it, really just that: a data structure. We also need data structures for another purpose, though, namely to represent our system's statements.

By an "algorithm", I mean *some kind of deterministic, orderly process* that can "process" data structures. It doesn't have to be a procedure that can be carried out by a Turing machine, however – we could imagine more powerful "machines", for example "machines" capable of carrying out so-called *supertasks* or *hypertasks* (i.e. tasks that require an infinite or even an uncountably infinite number of steps). What is essential here, though, is that we are talking about *orderly* processes. The process of determining if a statement is valid and the process of evaluating a statement must be orderly – that is an absolutely quintessential property of formal systems. In fact, I submit that for a linguistic system to be of any use to us, its formal grammar and evaluation function must not only be *orderly* processes – we must also be able to *understand* and *define* how they work. If we don't understand how they work, knowing the answer to a query wouldn't tell us anything useful about the World – it would be like feeding our statement to a black box and trying to make sense of its output when we have no idea of what it does.

Finally, we have the third core concept – semantics. The semantic schema defines our system: it defines how our model is to be constructed and what the data in it represents. Note, however, that it also needs to define the formal grammar and the evaluation function, i.e., how the language itself works.

## Many ways to create a formal system

If we accept this as our definition of what a "formal system" is, we must obviously conclude that there are many ways to construct a formal system. Many different types of formal systems have also been described in the literature – systems with different features, different syntaxes, and different expressive powers. The most commonly used formal language among philosophers and logicians is definitely *first-order predicate logic* (PL), but this is just a matter of popularity – PL is by no means the only imaginable formal language.

This calls for humility. For all I know, there could very well be an advanced civilization of squid-like aliens in some galaxy far, far away that have developed a language fundamentally different from our languages – a language based entirely on gestures and touch, say. Who knows what questions and problems occupy their philosophers and mathematicians?

## Anything goes?

The fact that there are many ways to imagine the construction of a formal system begs a lot of questions. What is "the right way" to design a formal system? What is "the right way" to do logic? Or are questions like these perhaps meaningless? Can we just *imagine* any formal system we want?

Well, we *are*, of course, indeed free to imagine whatever we want, but I submit that when we speak of "truth", we are not interested in imaginary fantasy worlds but in the *real* world – the notion of truth rests on the assumption that there is an objective Reality. We are thus actually *not* free to imagine whatever formal system we want – not if we are interested in *truth*.

## On objects, relations, and the principle of bivalence

As has long been noted, we humans tend to think in terms of *objects* and *relations* between objects. Conventional PL relies entirely on these two basic concepts, as does virtually every formal logic described in the literature. But is this "object/relation paradigm" the *only* way to describe the Reality?

I would say the answer is most definitely "no" – it is perfectly possible to imagine formal languages without quantifier-like mechanisms. Maybe our squid-like alien friends in the galaxy far, far away have a language that isn't object/relation-based – it is certainly possible, for all I know. I would argue, however, that the object/relation paradigm is *generic*, i.e., that it can be used to encode *any* kind of information. If this is the case, we can, without loss of generality, limit our investigation to object/relation-based languages.

A similar argument can be made when it comes to the principle of bivalence, i.e., the idea that logical statements are either true or false and nothing else. It is certainly possible to imagine formal systems where there are more than two possible truth-values – i.e., *many-valued logic*. However, it can be argued that many-valued logics can always be reduced to classical, two-valued logic. This is known as *Suszko's thesis* [3]. Some disagree with Suszko's argument and interpretation, e.g. Malinowski [4], but if we accept the information-theoretic view of language, I think it is hard to avoid the conclusion that indeed many-valued logic does not enable us to convey any information that could not be conveyed using bivalent truth-values. We can argue for this conclusion by noting that we, in practice, reasonably can only ever *absorb* a finite amount of information as answer to any query, and reasonably, any finite amount of information can be encoded using a finite string of bits.

In this paper, I am therefore going to assume that we indeed can describe the Reality using a bivalent logical language based on the object/relation paradigm. If that is not the case – well, that would invalidate pretty much the very idea of logic as we know it, so we have better believe that this is possible.

## The necessity of a self-contained system

What, then, about the formal grammar and the evaluation function? Must not they also be a part of our model? It seems to me that indeed they must. If we model the Reality as a bunch of objects and relations, then these objects and relations are all there is. If the evaluation function is not a relation, then what would it mean to say that a statement is true? It would be saying that this statement has a property that doesn't exist and doesn't correspond to anything real, and I would argue that this makes no sense. Nothing can "have" a non-existent, non-real property.

I thus contend that we need a formal system that is *self-contained*. The formal grammar and the evaluation function must be a part of the system's model. This, however, is where conventional PL and the standard axiomatic framework of mathematics – the ZF axioms – run into problems.

## Is it really that big a deal?

Some will now most definitely object, arguing that it isn't really a big deal if our system isn't self-contained. Tarski argued that we just have to accept that we cannot talk about the truth predicate of our language within the language itself. If we want to talk about it, we need another, stronger language – the meta-language. If we want to talk about the truth predicate of the meta-language, we need a meta-meta-language, and so on.

Thus, Tarski envisioned a hierarchy of languages, where each language can talk only about the languages on the lower levels. This is analogous to how ZF set theory avoids Russell's paradox: in ZF, we can think of the sets as organized in a hierarchy, where the sets on one level can have as members only sets from the lower levels.

This might seem a simple and natural solution to the problem, but I contend that it is far from a satisfactory answer.

Let us start with Tarski's hierarchy of languages and meta-languages. The moment we try to formalize this idea, we will run into problems. To formalize it, we need a formal language, but which? No language *in* the envisioned hierarchy can talk about the *entire* hierarchy. And if we use a language that isn't in the hierarchy – well, then our hierarchy is obviously incomplete.

Conventional set theory suffers from similar problems. If you ask a mathematician what a "relation", a "function", or a "model" is, you will likely be pointed to the standard definitions of these terms – definitions based on set theory. However, the set membership relation, $\in$, for example, is not a relation by these definitions. "Relations" defined by logical formulae ranging over all sets are also generally not proper relations by this definition. Moreover, by these definitions, our set theory itself has no model (this paradox I find particularly hard to brush away: if we set out to define the concept of a "model", and in doing so create a theory which does not have a model, it seems obvious to me that we have failed to give a complete and coherent definition of what a "model" is).

What about class theory? Well, it doesn't help. Just like ZF set theory, standard class theory (defined by the NBG axioms) is also well-founded, i.e., hierarchical – it just adds one more level to the hierarchy. We will end up with exactly the same problems; only this time phrased in terms of classes rather than sets.

This touches on another problem with well-founded set theory: the determination of the "size" of the hierarchy. No matter what universe of sets I imagine, I can always imagine a universe with even *more* sets by just extending the hierarchy. So how do we decide where to stop? No matter where we stop (and we *have* to stop somewhere), the choice of stopping there rather than somewhere else is going to seem arbitrary.

# Alternative approaches

The conventional mathematical foundational framework thus suffers from a number of problems – in fact, very *fundamental* problems, I contend – and there have been many attempts to find other approaches.

Kripke [5] proposed turning to three-valued logic. Propositions like the liar sentence are, according to this view, neither true nor false, but have a third truth-value: they are, as he called them, *ungrounded.* This might seem a simple and natural solution, however it runs into problems with the so-called "revenge" of the liar paradox. The approach allows us to talk about our *truth predicate*, but it does not allow us to talk about or define the notion of *groundedness*. It thus still suffers from a very fundamental form of undefinability.

Others have tried to resolve the problems with vicious circularity by proposing set theories that are not well-founded. The most popular of these non-well-founded set theories is definitely Quine's *New Foundations* (NF) [6]. Other systems have been proposed by Finsler [7], Boffa [8], and Forti & Honsell [9] (also known as Aczel's anti-foundation axiom). These alternative set theories have not gained much acceptance within the mathematical community, however, and though they are self-contained in some important respects (e.g., they include a set of all sets), they are all based on PL, which means they do not in themselves solve the problem with the undefinability of the truth predicate.

With this paper, I want to propose a new framework, which I like to call *Turing Verifiable Logic* (TVL), which approaches the problem in a different way.

## Turing Verifiable Logic

As the name suggests, Turing machines play a central role in the definition of TVL. The underlying idea is more generic than that, however, and can be extended to more powerful types of machines, for example machines capable of carrying out supertasks or hypertasks. TVL is still a very interesting case, though, as it seems to be the simplest possible system that is self-contained in this way.

Essentially, TVL postulates that by an "orderly process", we mean a process that can be carried out by a Turing machine, and that only objects and relations that can be created by such processes exist.

This is obviously not a new idea. The first to suggest that our Universe could be described as a Turing computation was Zuse [10]. More recently, similar thoughts have been voiced by authors like Schmidhuber [11] and Tegmark [12]. This general idea has been called *Zuse's thesis* (as suggested by Schmidhuber) or *the Computable Universe Hypothesis* [12]. However, we need to figure out how to deal with non-halting Turing machines and what to make of formal languages like PL, whose truth predicates can't be described as Turing computations. Both Zuse, Schmidhuber, and Tegmark have been vague on this point. Tegmark has voiced the concern that Gödel's incompleteness theorem might "torpedo" the whole idea [13], and similar questions have been raised by other authors, e.g. Vilenkin [14]. This is where I think this paper can offer new ideas.

# The TVL Ensemble

Let us start with the construction of the model. Our model, as we here choose to define it[1], is going to be a bunch of objects and relations – just like we're used to from PL. However, unlike in conventional PL, I am going to consider relations to be a special kind of object. In other words: every relation is also an object. There are however some objects that are not relations, and I will refer to them as *pure objects*.

Our model can thus be thought of as a bunch of objects, where each object is either a pure object or a relation. The difference between a pure object and a relation is that a relation expects a certain number of *arguments*. The number of arguments a particular relation $r$ expects is constant, and we will refer to this number as the *dimensionality* of $r$.

The model is going to contain an infinite number of pure objects and an infinite number of relations. Each object in our model will also be assigned a finite, unique *name*. The pure objects we will name using numerals, i.e., we name them 0, 1, 2, and so on. Indeed, as we shall see, we can think of them as numbers – and I will, therefore, often refer to them as "numbers" here.

The relations in our model are similarly assigned the names R0, R1, R2, and so on. These names we will refer to as "R-names". When we talk about specific relations, we will typically use a more convenient alias to denote the relation in question rather than its (probably quite long) R-name, but when I speak of an object's *name*, I always mean a numeral (for pure objects), or an R-name (for relations).

What we have in our model is thus the objects named 0, 1, 2, etc., and the relations named R0, R1, R2, etc. And that's it – that's all there is in our model!

Note that since I consider relations to be objects, relations can also be passed as input arguments to other relations. They can even be passed as input arguments to themselves, which is an interesting and quite useful feature. It is thus, for example, possible that the relation R17 might hold when given the objects 5 and R0 as input, in which case we say that the relational expression "R17(5, R0)" is true. Similarly, it could be the case that the relation R0 holds when given *itself* as input, i.e., that "R0(R0)" is true, and so on.

An example of a situation where it is convenient to be able to use relations as objects is when we want to talk about the dimensionality of relations. We would like to have a relation *DIM*($o$, $n$), which holds true if and only if the object $o$ has a dimensionality of $n$ (for convenience, let us say that pure objects have a dimensionality of 0). Such a relation indeed exists in our model, and it assigns a dimensionality to every single object in the Ensemble.

---

[1] There are, of course, many ways to encode the same information, and consequently there are many effectively equivalent ways to describe the same thing. Here we have to choose a particular way of defining our model and our language, but the underlying idea could also be "realized" in many other ways.

But *DIM* is also itself a relation of dimensionality 2, which means that it is the case that *DIM*(*DIM*, 2).

## The construction

To construct the TVL Ensemble, we need a programming language, i.e., a language for defining Turing algorithms. This language, which we will henceforth denote $\mathcal{L}$, should be Turing complete – that is, it should allow us to define *any* procedure which can be carried out by a Turing machine. $\mathcal{L}$-programs may also take input arguments, and for simplicity, let us decide that every $\mathcal{L}$-program should begin with a *signature declaration*, which defines how many input arguments this program expects (the number can be 0). The input arguments to our algorithms are always strings. This is all we need to say about $\mathcal{L}$ for now – as long as $\mathcal{L}$ has these properties, we are free to design $\mathcal{L}$ pretty much any way we want.

Consider now the collection of all programs written in the programming language $\mathcal{L}$ which take one or more input arguments and order these programs alphabetically. This gives us an infinite sequence of strings, $A_0$, $A_1$, $A_2$, and so on, defining Turing algorithms. The idea is now that the $i^{th}$ relation, R$i$, should hold when given the objects $x_1, x_2, ..., x_n$ as input arguments if and only if the $i^{th}$ algorithm, $A_i$, takes $n$ input arguments and it *halts* if it is given the *names* of the objects $x_1, x_2, ..., x_n$ as input. So, for example, R17(5, R0) is true if and only if the algorithm $A_{17}$ takes 2 arguments and it halts if it is given the strings "5" and "R0" as input.

That's basically it: this is how we construct the TVL Ensemble.

Note that we are not interested in any kind of output values from the algorithms – we are only interested in whether or not they halt. This is a very crucial point.

## The language

Our description system also needs a formal language so that we can talk about our Ensemble of objects, and we can actually use our programming language, $\mathcal{L}$, as this formal language.

By a *closed TVL-statement*, we simply mean a string that defines an $\mathcal{L}$-program which takes no input arguments, and we say that this statement is *true* if and only if the algorithm in question *halts*.

We can also generalize this idea and allow for *open* TVL-statements. Open TVL-statements work just like the relations in our Ensemble: by an open TVL-statement P, we mean a string that defines an $\mathcal{L}$-program, taking some number of input arguments $n$. This P is *true* for the input arguments $x_1, x_2, ..., x_n$ if and only if the program defined by P halts when it is given the *names* of the objects $x_1, x_2, ..., x_n$ as input.

Note that this means that the behavior of every single relation in the Ensemble can be described using an expression formulated in the TVL language (i.e., $\mathcal{L}$). That is: for every relation R$i$ in the Ensemble, there is an open TVL-expression that holds for precisely those input arguments which make R$i$ true. All relations in the TVL universe are thus definable – a very nice property.

## Consistency

The construction of the TVL Ensemble and the assignment of truth-values to TVL-statements is quite straightforward: basically, it all comes down to determining which Turing machines halt and which do not. If we had access to a *TVL-oracle* – a magical "machine" capable of determining instantly whether or not a particular $\mathcal{L}$-program will halt if given certain input arguments – we could use this oracle to figure out which relational expressions are true and which are not, and which TVL-statements are true and which are not.

It is also well-known that Turing algorithms can be described using number theory. Any proposition describing a halting problem, i.e. any proposition of the form "The $\mathcal{L}$-program P halts if it is given input arguments so-and-so", can be mechanically translated into a proposition about numbers – that is: a proposition of conventional (i.e., PL-based) number theory. Having access to a "number-theory-oracle" would thus allow us to predict all answers provided by a TVL-oracle. In other words: if number theory is consistent, then so is the construction of the TVL universe and the TVL truth predicate.

We can thus be quite confident that the construction is consistent (indeed, it is hard to get more confident than that when it comes to consistency proofs).

## A conservative but rich ontology

The TVL Ensemble will obviously contain an infinite number of objects. However, it will "only" contain a *countably* infinite number of objects (i.e., a set theorist would say that the cardinality of our model is $\aleph_0$). This is a lot "fewer" objects than the ZF universe supposedly contains. The TVL Ensemble only contains objects which can be constructed using very simple means, and thus presents a very conservative mathematical ontology.

Nevertheless, the TVL universe contains many interesting relations. For example, it is easy to see that in the TVL Ensemble, there will be a relation NUMBER(*o*) which holds if and only if *o* is a pure object – i.e., what we mean by a "number" in the TVL world. Defining this relation is easy: the corresponding algorithm is going to be given the name of *o* as input, and all it needs to do is to check if this name starts with an "R". If it does not start with an "R", *o* must be a pure object, i.e., a "number", so we let our algorithm halt; if it *does* start with an "R", however, we let our algorithm step into an infinite loop. This algorithm will thus halt if the input is the name of a "number" and loop forever if it is the name of a relation. The corresponding relation NUMBER(*o*) will thus be true if *o* is a number and false if *o* is a relation – which is precisely what we wanted.

(Note how we, in this example, deliberately let our algorithm step into an infinite loop. This may seem a bit odd at first, but we will be doing this a lot since the *only* way for a relation or a TVL-statement to evaluate to false is if the corresponding algorithm does not halt.)

Similarly, it is easy to see that in the TVL Ensemble, there will be relations corresponding to basic arithmetic operations (addition, multiplication, exponentiation, and so on), acting on the pure objects. That is: there is, for example, a relation ADD(*a*, *b*, *c*) which holds true if and only if *a*, *b* and *c* are "numbers" and $a + b = c$. This also justifies thinking of the pure objects in our Ensemble as "numbers".

We can obviously also use these "numbers" to represent strings, using Gödel-numbering techniques. We can thus discuss strings, statements, and grammar.

It is also clear that in the Ensemble, there will be a binary relation, DIM, which works as discussed before (i.e., determines the dimensionality of a given object). The construction is a bit more complicated, so let us say something more about it. We want this relation DIM(*o*, *n*) to hold if and only if *n* is a number corresponding to the dimensionality of the object *o*. We thus construct an algorithm that begins by checking that the second input argument is a number. If not, we can immediately step into an infinite loop (i.e., make the relation evaluate to false since the second argument is of an invalid type). If *n* is a number, on the other hand, we need to examine the object *o*. If *o* is also a number, we know that its dimension is 0. If it is a relation, we need to compute which algorithm it corresponds to. This should be straightforward: we would have the R-name of *o* as our input, i.e., we would have its name, of the form R*i*, and from there, we can (with a bit of effort) find the code of the corresponding algorithm $A_i$. Since we've decided that our programming language $\mathscr{L}$ requires all programs to begin by declaring the number of arguments they take, we can then analyze the code of $A_i$ to extract this number, which is also the dimensionality of the relation R*i*. Thus we can determine the dimensionality of *o*, regardless of whether *o* is a number or a relation, and we can compare this number to *n* and decide whether we should halt or step into an infinite loop – and that is all we need to construct the DIM-relation.

So while the TVL Ensemble represents a very conservative mathematical ontology (compared to conventional set theory), it is nonetheless clearly a very rich world.

## Defining the truth predicate of TVL using TVL

The most remarkable property of TVL, though, is that we can use the TVL language (i.e., $\mathscr{L}$) to define the truth predicate of the TVL language itself. That is: there is a definable

TVL-relation *TRUE(s)* which holds if and only if *s* is a string (i.e., a Gödel number) defining a true, closed TVL-statement.

The construction is simple. The corresponding algorithm $A_{TRUE}$ begins by checking if the input argument *s* is a Gödel number corresponding to a valid, closed TVL-statement. If not, we can step into an infinite loop immediately (which means that *TRUE(s)* will evaluate to false for input arguments of invalid format). If *s* indeed is a valid statement, i.e., if *s* corresponds to a definition of an algorithm *A* expressed in the language $\mathcal{L}$, we instead just execute *A*. This is not a problem: it is well known that Turing algorithms can execute other Turing algorithms provided their code as input. This is achieved by creating a *virtual* Turing machine within the current Turing machine. We then just execute this algorithm *A* until it halts – if it ever does. If it does halt, we know that the statement *s* is a true TVL-statement, so we also let our algorithm ($A_{TRUE}$) halt. If *A* never halts – well, then we will just continue indefinitely to execute *A* in our virtual Turing machine, and in this case, $A_{TRUE}$ will obviously never come to a halt either. That is: $A_{TRUE}$ will halt if and only if *s* is a valid program that halts – and this is *exactly* how we need the algorithm corresponding to *TRUE(s)* to work.

Note that we don't need our algorithm $A_{TRUE}$ to predict within a finite amount of time if the given algorithm, defined by *s*, will halt or not (that would be impossible, as Turing demonstrated). We just need to design $A_{TRUE}$ so that it *co-halts* with the algorithm defined by *s*, and this is perfectly possible and does not violate Turing's results. We can thus use the TVL language to define the truth predicate of the TVL language itself – and this truth predicate *TRUE(s)* also exists in our Ensemble, as one of the relations, with some R-number.

## What about the liar paradox?

So what about the liar paradox? It does not kick in, because it is simply not possible to construct a statement in TVL which says "This TVL-statement is false". As we have seen, we *can* define an algorithm $A_{TRUE}$ which evaluates to true if and only if it is given a true TVL-statement as input, but we can *not* define an algorithm $A_{FALSE}$ which does the opposite. We can write a Turing algorithm that emulates and co-halts with a given input algorithm, but we cannot write an algorithm that halts if and only if the input algorithm does *not* halt (that would require us to solve the halting problem, which we know is impossible).

Our Ensemble thus has a property that identifies all true statements, but there is no property that identifies all false statements. Rather, falseness must be understood as the *lack* of truthness. This also means that, in general, TVL-statements cannot be negated. Indeed, TVL is altogether a weaker system than conventional number theory. As we noted before, we could predict the output of a TVL-oracle to *any* query if we had access to a number-theory-oracle, but the reverse is not true.

# An alternative formulation

So far, we have used $\mathcal{L}$, a programming language, as our formal language. TVL-statements are thus definitions of programs, and whether or not a TVL-statement is true depends on whether or not the corresponding program halts. This is obviously very different from PL. In $\mathcal{L}$, we have no quantifiers, logical operators, or relational expressions in the usual sense (recall that I mentioned earlier that it is possible to conceive of languages that do not have any quantifier-like mechanisms – here we have a good example of such a language).

It is, however, also possible to define an alternative formal language for talking about the TVL Ensemble, a language that has the same expressive power as $\mathcal{L}$, but which is more similar to PL in structure (though with some additional restrictions). Let us have a look at this alternative formulation, as it will allow us to see the differences between TVL and PL more clearly.

This alternative, "PL-like" language I will here denote $\mathcal{L}^*$, and it is defined by the following recursive rules:

**Terms:** A *term* refers to an object (i.e., either a pure object or a relation). It may be the *name* of a specific object or a *variable*.

**Atomic expressions:** Atomic expressions are primitive expressions about terms. The valid atomic expressions we have are statements of equivalence ($t_1 = t_2$), statements of non-equivalence ($t_1 \neq t_2$), and relational expressions ( $t_0(t_1, t_2, ..., t_n)$ ). By definition, let us say that a relational expression evaluates to *false* if the number of arguments provided doesn't match the dimensionality of the relation in the expression.

**Arithmetical expressions:** For convenience, let us also add simple arithmetical statements about terms as "built-ins" of our language, i.e., expressions like "$2 + 2 = 4$", "$x + 5 = y$", "$17 \cdot z > 36$", etc. If a term that occurs in an arithmetical expression refers to a relation rather than a number, we can say that the whole expression evaluates to false.

**Logical connectives:** We can use the logical connectives $\wedge$ ("and") and $\vee$ ("or") between any two $\mathcal{L}^*$-expressions. We can, however, *not* use other PL connectives like $\neg$ ("not"), $\rightarrow$ ("implies"), or $\leftrightarrow$ ("if and only if"). This is because in TVL, we cannot use negations without restrictions. The $\rightarrow$ and $\leftrightarrow$ must also be forbidden because they contain "hidden" negations (if they were to be allowed, they could be used to negate arbitrary expressions).

**Quantifiers:** We can use the existential quantifier ($\exists$) as usual. That is: given a $\mathcal{L}^*$-expression $P$ with a free variable $x$, the statement "$\exists x : P(x)$" is also a valid expression. We can *not* use the universal quantifier ($\forall$) as usual, though – but we can use a *bounded* version of it, which we can write: "$\forall x \leq b : P(x)$", meaning that all of the objects 0, 1, 2, ..., $b$ and R0, R1, R2, ..., R$b$ satisfy $P$. If the term $b$ doesn't refer to a number, we define that this $\forall$-expression evaluates to false by default.

It is fairly straightforward to prove that ℒ and ℒ* have exactly the same strength: any ℒ*-statement can be mechanically translated to an ℒ-statement, and vice versa (see the appendix for a proof sketch).

The principal difference between TVL and PL is thus basically that in TVL, we only have the equivalence of *bounded* ∀-quantifiers and the logical connectives ∧ and ∨.

# The limits of TVL

TVL might seem like a very *weak* language with all these restrictions. I contend, however, that it is actually not as weak as it may seem at first glance. In fact, I would argue that if we had access to a TVL-oracle, we could settle most of the questions that occupy today's mathematicians.

## To believe or not to believe

Let us start by considering a simple proposition, which at first might seem impossible to formulate using TVL, namely the Goldbach conjecture, G, which asserts that every even integer greater than 2 can be written as the sum of two primes. It may seem that G cannot be expressed in TVL since it has an unbounded universal quantifier. However, we *can* express the *negation* of G in TVL. To do that, we just need to write an algorithm that loops over all the integers looking for a *counterexample* and halts if it finds one.

If we had access to a TVL-oracle, we could thus just feed it this negative form of the proposition. If the oracle gave us a positive answer (i.e. if the algorithm finds a counterexample and halts), we would know that Goldbach was wrong; if the oracle gave a negative answer (i.e. if the algorithm never comes to a halt), we would know that there is no counterexample, and thus that Goldbach was right.

So in a way, we, so to speak, have "one layer of negations", since our TVL-statements are ultimately either true or false, and we can either believe in a certain statement or not. Knowing that a TVL-statement is *false* also tells us something about the TVL Ensemble.

## Delimiting universal quantifiers

Now let us consider another example, namely Euclid's theorem, E = "For every prime $p$, there is a larger prime $q$". In this case, the counterexample technique we used for the Goldbach conjecture does not seem to work, but there is another trick we can use: we can try delimiting our universal quantifiers.

Conventional wisdom tells us that for every prime $p$, there should be another prime $q$, which is larger than $p$ but less than or equal to $p! + 1$. There is thus a *computable limit* (defined by the factorial function plus 1), such that we don't need to look further than this limit before we should find another prime. Now consider the statement NE' = "There exists a prime $p$ such that no prime exists between $p + 1$ and $p! + 1$". This statement *can* be expressed in TVL, and if we could feed it to a TVL-oracle, we would expect it to return false – and if we find that NE' is indeed false, we would also effectively have demonstrated that Euclid was right (actually, NE' being false even tells us a bit more than the statement E, since it also tells us something about how far apart the primes are from each other).

## Stronger systems

Using "tricks" like these, I would expect that a mathematician with access to a TVL-oracle would be able to settle most open questions in contemporary mathematics. Moreover, even though TVL is a weaker system than conventional PL-based number theory, the idea behind the construction of the TVL world can be generalized to yield *stronger* systems that are self-contained in the same way.

Suppose, for example, that instead of regular Turing machines, we use *oracle Turing machines* – that is: Turing machines equipped with an *oracle function*, which is capable of calculating some function which cannot be calculated by regular Turing machines. Thus we can construct different flavors of *Oracle Turing Verifiable Logic* (OTVL). We could, for example, imagine an oracle Turing machine with an oracle function that takes a string as input and returns true if and only if this string represents a true proposition of conventional number theory (that is: a proposition about numbers, expressed using *classical* PL). We then similarly construct a programming language for this new "machine", call it $\mathscr{L}^O$, which is just like $\mathscr{L}$, except it also has a command for invoking the oracle function. Then we construct our OTVL system the same way as before.

The new system is clearly at least as strong as conventional number theory (given access to an oracle for our new system, which can tell us which OTVL-statements are true and not, we could obviously predict all answers from a number-theory-oracle, by invoking the oracle function). But we also get a logical language that is self-contained, because just like regular Turing machines, oracle Turing machines can emulate other oracle Turing machines of the same type (i.e., with the same oracle function).

Self-containment "TVL-style" thus doesn't necessarily require a system weaker than number theory. Indeed, it would seem we could take *any* model that we may envision and imagine a self-contained language constructed from it – we can always, so to speak, add an "OTVL-layer" on top of it, with the truth predicate of the original system as the oracle function.

# What is an "orderly process"?

This brings me back to the question of what exactly substitutes an "orderly process". This is a very fundamental question – indeed, I would argue that the search for the "right" foundational framework for mathematics *is* the search for the "right" answer to this question.

In the TVL world, an "orderly process" is a process that can be carried out by a Turing machine (possibly one that never halts). To evaluate the PL truth predicate of number theory, however, we would need a much more powerful machine: we would need a machine that can carry out a supertask – a task that requires *more*[2] than an infinite number of steps. In set-theoretic terms, a Turing machine can be described as a machine capable of carrying out $\omega$ operations, whereas a "machine" for calculating the PL truth predicate for number theory would require $\omega^\omega$ operations. A machine for evaluating the truth predicate of the OTVL-system we described before would require an even more powerful machine, one that can carry out $\omega^\omega \cdot \omega$ steps.

We could go on, of course, imagining even more powerful machines. However, as we noted before, no matter what universe of sets you imagine, you can always imagine another, larger universe with more sets and more ordinals in it. Thus, no matter what "machine" you imagine, it is always going to be possible to imagine a more powerful one – one that can carry out even more steps, a machine that can "talk" about the first machine.

TVL can thus be understood as the consequence of *rejecting the idea that an "orderly process" can carry out "more" than an infinite number of steps*.

# A potential Theory Of Everything?

Is it conceivable that the TVL Ensemble, when interpreted by some semantic schema, could describe the Reality completely – that is: could TVL be a *Theory Of Everything* (TOE)?

Well, I think it is possible. Indeed, considering all the problems the conventional framework suffers from, I am inclined to think that TVL is at least a more promising candidate for a TOE than ZF set theory is.

At any rate, I think we can safely dismiss Tegmark's concern about Gödel's incompleteness theorem being a potential deal-breaker for the Computable Universe Hypothesis. Gödelian incompleteness is a problem in PL and PL-like languages, but PL is not the only conceivable formal language, nor is it obvious that PL is, so to speak, "the language of nature" (or "God's language", if you will). If the TVL Ensemble is a way to describe the Reality completely, then the truth predicate of PL would simply not exist – but the truth predicate of the TVL language *would* exist, and as we have seen, this language would be nicely self-contained.

---

[2] In the ordinal sense, not the cardinal sense.

# Conclusion

Virtually all mathematical research today is based on ZF set theory and PL, but there are many fundamental problems with this framework, and mathematicians still debate how to interpret the discoveries of Russell, Gödel, and Tarski.

In this paper, I have proposed a different system, TVL, with a formal language that can talk about and define its own truth predicate. Such a formal language has not been described in the literature before, as far as I know.

TVL, however, requires us to accept some restrictions compared to what we are used to from conventional PL. On the other hand, *all* formal languages, including PL, have limits when it comes to what they allow us to express. If our language doesn't allow us to formulate the equivalent of all types of PL expressions, then that is obviously a limitation of our language; but if it *does* allow us to formulate all such expressions, it cannot talk about its own truth predicate, and that would *also* be a limitation.

We are thus caught between a rock and a hard place: we need to relinquish one of these seemingly very desirable properties, for we can't have both – but which one should we give up on? So yes, TVL indeed has some features that don't seem ideal, but so does the conventional framework.

TVL can be understood as the consequence of rejecting the idea that an "orderly process" can carry out "more" than an infinite number of steps, which I would argue is a rather natural limitation. We can, of course, *imagine* "processes" and "structures" that continue beyond infinity, and this is also exactly what we do when we use PL and ZF set theory. It doesn't seem to lead to any outright contradictions (in the PL sense), but it opens quite a worm can of problems and paradoxes.

If we believe that there is some kind of mathematical Reality (as I do), i.e., if we are *platonists*, then the ultimate question is, of course, not what kind of structures we can *imagine* but what kind of structures correspond to something *real*. Considering all the problems we have discussed, I personally strongly doubt that the extremely inflated ontology postulated by conventional set theory corresponds to something real, and I would argue that TVL is a more promising candidate for a TOE.

At any rate, TVL is an interesting example of what an alternative mathematical ontology and formal language might look like. It demonstrates that it is possible to construct a formal language that can define its own truth predicate, which is an interesting result in itself.

# APPENDIX

**Proof (sketch) that it is possible to mechanically translate any statement formulated in $\mathcal{L}$ into an equivalent statement in $\mathcal{L}$\*, and vice versa:**

(The proof is fairly straightforward and only relies on well-known techniques, so only a sketch is given).

**From $\mathcal{L}$\* to $\mathcal{L}$:** Obviously, every atomic $\mathcal{L}$\*-expression, P\*, can be mechanically translated to a corresponding $\mathcal{L}$-expression, P, where P defines an algorithm that takes one input argument for every free variable in the expression P\*, and halts if and only if P\* is true for the given arguments. Further, if P\* and Q\* are two $\mathcal{L}$\*-expressions (possibly with free variables), and they can be mechanically translated to $\mathcal{L}$-expressions P and Q, then we can also mechanically translate the $\mathcal{L}$\*-expressions "P\*$\wedge$Q\*", "P\*$\vee$Q\*", "$\exists x : P^*(x)$", and "$\forall x \leq b : P^*(x)$" into equivalent $\mathcal{L}$-expressions. For "P\*$\wedge$Q\*" we create a new algorithm that just executes P and Q in sequence (thus halting if both P and Q halt); for "P\*$\vee$Q\*" we create an algorithm that executes P and Q in *parallel*, in two internal virtual Turing machines (VTMs), halting if at any point either P or Q halts; for the bounded $\forall$ quantifier we execute P(x) for each value of x in the range (which is finite) in a sequential loop, halting if P(x) halts for each and every value for x in the range; for the $\exists$ quantifier, we create a growing pool of VTMs, each running P(x) for some value of x, and continue to expand the pool with VTMs for more values of x, while also executing the VTMs in the pool, and we halt if any of the VTMs in the pool halts at some point. By recursively applying these translations, we can translate any $\mathcal{L}$\*-statement into an equivalent $\mathcal{L}$-statement.

**From $\mathcal{L}$ to $\mathcal{L}$\*:** Design a $\mathcal{L}$\*-expression EXEC-SEQ($p$, $x$) which is true if and only if $p$ is a string (i.e. Gödel number) defining an $\mathcal{L}$-algorithm without input arguments, and $x$ is a string (i.e. Gödel number) which demonstrates that the algorithm $p$ halts. Design it so that to satisfy this formula, $x$ must be a string divided into "frames", separated by some special delimiter character, where each subsequent "frame" defines a subsequent step of the execution of $p$. Each frame should hold information about what the content of the memory tape is and which the current program instruction is at the corresponding stage in the execution. The formula should also check that the final frame corresponds to the program having reached the instruction to halt. Checking if a given $p$ and $x$ satisfy the necessary internal conditions does not require any unbounded loops, and it should thus be possible to check this with a $\mathcal{L}$\*-expression. Using this, we can translate any closed $\mathcal{L}$-statement P into an equivalent $\mathcal{L}$\*-statement P\* by using an expression of the form "$\exists x : $ EXEC-SEQ($p$, $x$)", where $p$ is simply hard-coded to the Gödel number of P. With minor changes, we could expand this to also cover open statements.